\documentclass[12pt]{article}

\usepackage{a4wide}
\usepackage{amsmath}
\usepackage{amsthm}
\usepackage{amssymb}
\usepackage[english]{babel}
\usepackage{graphicx}

\newcommand{\Z}{\mathbb{Z}}

\newcommand{\gr}{\textup{gr}} 

\DeclareMathAlphabet{\mathcalligra}{T1}{calligra}{m}{n}

\DeclareFontFamily{U}{mathx}{\hyphenchar\font45}
\DeclareFontShape{U}{mathx}{m}{n}{
      <5> <6> <7> <8> <9> <10>
      <10.95> <12> <14.4> <17.28> <20.74> <24.88>
      mathx10
      }{}
\DeclareSymbolFont{mathx}{U}{mathx}{m}{n}
\DeclareFontSubstitution{U}{mathx}{m}{n}
\DeclareMathAccent{\widecheck}{0}{mathx}{"71}

\newtheorem{theo}{Theorem}[section]
\newtheorem{lemma}[theo]{Lemma}
\newtheorem{prop}[theo]{Proposition}

\newtheorem{rem}[theo]{Remark}
\newtheorem{ex}[theo]{Example}

\newtheorem{defi}[theo]{Definition}

\begin{document}
\title{Spectral sequences of a Morse shelling}
\author{Jean-Yves Welschinger}
\maketitle

\begin{abstract} 
\vspace{0.5cm}

We recently introduced a notion of tilings of geometric realizations of finite relative simplicial complexes and related those tilings to the discrete Morse theory of R. Forman, especially when they have the property of being shellable, a property shared by the classical shellable complexes. We now observe that every such tiling supports a quiver which is acyclic precisely when the tiling is shellable and then, that every shelling induces two spectral sequences which converge to the relative (co)homology of the complex. Their first pages are free modules over the critical tiles of the tiling. \\


\end{abstract}

\section{Introduction}

We recently introduced \cite{SaWel1, SaWel2} a notion of tiling of the geometric realization of a finite simplicial complex. It is a partition by tiles, where a tile is a maximal simplex deprived of several of its codimension one faces together with possibly a unique face of higher codimension, see \S \ref{sectiles}. When the tiles can be totally ordered in a way defining a filtration by subcomplexes, the tiling is said to be shellable, or rather Morse shellable, since this notion extends the classical notion of shellability \cite{BruMan,AdipBen2,Stanbook,Z,HachZ,BenZ} and is related to discrete Morse theory \cite{For,SaWel2,Koz2}. In particular, the boundary of any convex simplicial polytope is shellable by \cite{BruMan}, but every closed triangulated surface is Morse shellable as well \cite{SaWel2} and 
in fact every finite simplicial complex becomes Morse shellable after finitely many stellar subdivisions at maximal simplices or after a single barycentric subdivision by \cite{Wel2}. Not all tilings are shellable though, for any product of a sphere with a torus of positive dimension carries Morse tiled triangulations which cannot be shelled, see \cite{Wel1}. 

We now observe that every Morse tiling supports a quiver which encodes its shellability, see \S \ref{subsecshellability}. 

\begin{theo}
\label{theointro1}
A Morse tiling is shellable if and only if its quiver does not contain any oriented cycle.  
\end{theo}

We then observe that every such shelling recovers the homology and cohomology of the complex, via two spectral sequences, see \S \ref{secspectral}.

\begin{theo}
\label{theointro2}
Any Morse shelling on a finite relative simplicial complex induces two spectral sequences which converge to its relative homology and cohomology respectively and whose first pages are free graded modules over the critical tiles. 
\end{theo}

A critical Morse tile is a simplex which has been deprived of several codimension one faces together with its remaining face of maximal codimension, see \S \ref{sectiles}. This terminology originates from its relation with the discrete Morse theory of Robin Forman \cite{For,Koz2}, for any Morse shelling encodes a class of compatible discrete Morse functions whose critical points are in one-to-one correspondance with the critical tiles of the shelling, preserving the index, see \cite{SaWel2} and Theorem \ref{theoMorse}. As in the case of discrete Morse complexes, the chain and cochain complexes appearing in the spectral sequences given by Theorem \ref{theointro2} have much lower ranks than the simplicial ones, see Remark \ref{remfinale}.

Finally, the quiver of a Morse tiled relative simplicial complex $S=K \setminus L$ makes it possible to weaken the shellabiity condition into some partial shellability which implies the Morse shellability of its $q$-skeleton, see Theorem \ref{theopartialshel}, and provides spectral sequences which converge to the relative (co)homology of the pair $(K,L)$ in degrees less than $q$, see Theorem \ref{theospectral}.\\

We compute in \S \ref{sectiles} the relative (co)homology of Morse tiles, showing that it does not vanish only for critical tiles. We introduce the quivers of Morse tilings in \S \ref{secquiver} and prove Theorem \ref{theointro1}. We also define the notion of partial shellability and prove Theorem \ref{theopartialshel}. We finally prove Theorem \ref{theointro2} in \S \ref{secspectral} and relate the filtration of a Morse shelling with discrete Morse theory in \S \ref{secMorse}.\\

\textbf{Acknowledgement:}
This work was partially supported by the ANR project MICROLOCAL (ANR-15CE40-0007-01).

\section{Relative homology of a Morse tile}
\label{sectiles}

Let $n$ be a non-negative integer. We recall that an $n$-simplex is the convex hull of $n+1$ points affinely independent in some real affine space or rather, abstractly, just a set of cardinality $n+1$ whose elements are vertices, see \cite{Munk,Koz}. A face of a simplex is the convex hull of a subset of its vertices and its dimension is the dimension of the affine space it spans. An open simplex (resp. face) is the relative interior of a simplex (resp. face) in its supporting affine space. The following definition has been given in \cite{SaWel2}.

\begin{defi}[Definition $2.4$ of \cite{SaWel2, Wel1}]
\label{defMorsetile}
A Morse tile $T$ of dimension $n$ and order $k \in \{0, \dots , n+1\}$ is an $n$-simplex $\overline{T}$ deprived of $k$ of its codimension one faces together with one possibly empty face $\mu$ of higher codimension. It is critical of index $k$ iff $\dim \mu = k-1$. The simplex $\overline{T}$ is called the underlying simplex while $\mu$ is called its Morse face. 
\end{defi}

\begin{figure}[h]
   \begin{center}
    \includegraphics[scale=1]{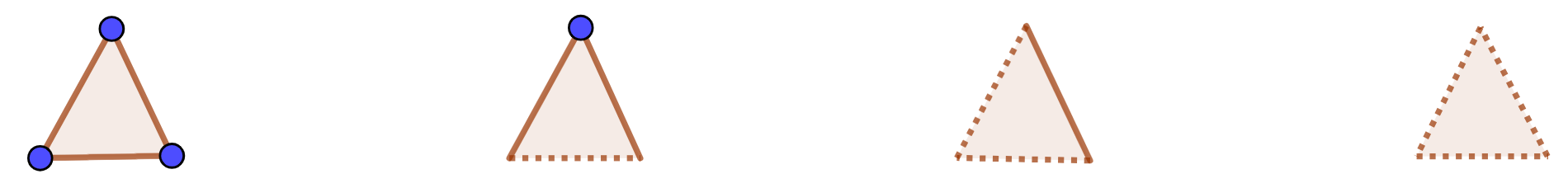}
    \caption{The basic tiles in dimension two.}
    \label{figbasictiles}
      \end{center}
 \end{figure}

A tile with empty Morse face is said to be {\it basic} and a basic tile of order $k >0$ contains a unique open face of dimension $k-1$, its restriction set \cite{Stanbook,Z}, and none of lower dimension, see Figure \ref{figbasictiles}. Critical tiles are thus the ones for which this peculiar face has been removed, see \cite{SaWel2,Wel2}. Note that critical tiles of maximal index are open simplices. The tiles which are not critical are said to be {\it regular}.

\begin{ex}
A tetrahedron deprived of one triangle and one edge is a regular Morse tile of order one and dimension three, whereas a tetrahedron deprived of two triangles and one edge is critical of index two. 
\end{ex}

Recall that a finite simplicial complex is a finite collection of simplices which contains all faces of its elements and such that the intersection of any two of them is a common face, possibly empty. Such a simplicial complex $K$ inherits a topology whose homology and cohomology can be computed out of the simplicial chain and cochain complexes
$C_* (K ; \Z) = \oplus_{\sigma \in K} \Z_\sigma$ and $C^* (K ; \Z) = \oplus_{\sigma \in K} \hom (\Z_\sigma , \Z)$. 
These complexes are graded by the dimension of the simplices and for every $\sigma \in K$, $\Z_\sigma$ denotes the infinite cyclic group whose generators are the two orientations on $\sigma$, that is the kernel of the augmentation morphism $a \sigma + b \overline{\sigma} \in \Z \sigma \oplus \Z \overline{\sigma} \mapsto a + b \in \Z$, see \cite{Munk}. 

\begin{defi}
\label{defrelhom}
The relative homology (resp. cohomology) of a Morse tile $T$ with coefficients in the Abelian group $G$ is the homology $H_* (\overline{T} , \overline{T} \setminus T ; G)$ (resp. $H^* (\overline{T} , \overline{T} \setminus T ; G)$) of the chain complex $C_* (\overline{T} ; \Z) \otimes G / C_* (\overline{T} \setminus T ; \Z) \otimes G$ (resp. $\ker \big( \hom (C^* (\overline{T} ; \Z) ,G) \to \hom (C^* (\overline{T} \setminus T ; \Z) , G) \big)$).
\end{defi}

\begin{prop}
\label{proprelhom}
The relative homology and cohomology of regular Morse tiles vanish whereas they are isomorphic to the coefficient group and concentrated in degree $k$ in the case of critical tiles of index $k$.
\end{prop}

\begin{proof}
Let $T$ be a Morse tile of dimension $n$ and order $k$, with underlying simplex $\overline{T}$ and Morse face $\mu$. Let $\sigma$ be the union of the codimension one faces that have been removed from $\overline{T}$, so that $T = \overline{T} \setminus (\sigma \cup \mu)$. Let $\theta$ be the convex hull of the $k$ vertices opposite to those faces, so that $\theta \subset \mu$. The tile $T$ is then regular iff $\theta \neq \mu$.

If $T$ is regular, then $\mu$ contains a vertex $v$ that does not belong to $\theta$, so that $\sigma \cup \mu$ is a cone with apex $v$. The inclusions $v \to \sigma \cup \mu \stackrel{i}{\to} \overline{T}$ then induce isomorphisms in homology and cohomology, these being concentrated in degree zero, see \cite{Munk}. The short exact sequences
\begin{equation}
\label{eqn1}
0 \to C_* (\sigma \cup \mu ; G) \stackrel{i_*}{\to}  C_* (\overline{T} ; G) \to C_* (\overline{T} , \sigma \cup \mu ; G) \to 0
\end{equation}
and
\begin{equation}
\label{eqn2}
0 \to C^* (\overline{T} , \sigma \cup \mu ; G)  \to  C^* (\overline{T} ; G) \stackrel{i^*}{\to} C^* (\sigma \cup \mu ; G) \to 0
\end{equation}
induce long exact sequences in homology with coefficients in the Abelian group $G$, from which we deduce that $H_* (\overline{T} , \sigma \cup \mu ; G) $ and $H^* (\overline{T} , \sigma \cup \mu ; G) $ vanish. 

If $T$ is critical, the result is well known for $k \in \{ 0, n\}$. When $k=0$, the relative (co)homology of $T$ is the (co)homology of a simplex while when $k=n$, it is the (co)homology of a simplex relative to its boundary. Otherwise, we deduce by excision that the relative (co)homology of the pair $(\sigma \cup \mu , \sigma)$ is isomorphic to the (co)homology of $\mu$ relative to its boundary, so that it is likewise isomorphic to the coefficient group $G$ and concentrated in degree $k-1$.  From the long exact sequence of this pair $(\sigma \cup \mu , \sigma)$, we deduce that the reduced homology and cohomology groups $\widetilde{H}_* (\sigma \cup \mu ; G)$ and $\widetilde{H}^* (\sigma \cup \mu ; G)$ get isomorphic to $H_* (\sigma \cup \mu , \sigma ; G)$ and $H^* (\sigma \cup \mu , \sigma ; G)$, so that they are isomorphic to $G$ and concentrated in degree $k-1$ as well.  The result then again follows from (\ref{eqn1}) and (\ref{eqn2}), for the reduced homology and cohomology of $\overline{T}$ vanish.
\end{proof}

\section{Quiver and shellability of a tiling}
\label{secquiver}

\subsection{Acyclic quivers}
\label{subsecacyclic}

\begin{defi}
\label{defquiver}
A quiver is a quadruple $Q=(Q_V ,  Q_A , s , t)$, where $Q_V$, $Q_A$ are finite sets whose elements are vertices and arrows, and where $s,t : Q_A \to Q_V$ are maps called source and target.
\end{defi}

A quiver is thus a finite directed graph which may contain several arrows between two vertices or even a loop arrow from a vertex to itself. 

\begin{ex}
\label{exquiver}
Let $(V , \leq)$ be a finite partially ordered set. It supports a quiver, its comparability graph, having $ V$ as set of vertices and containing an arrow for each pair $(i,j) \in V^2$ such that $i \leq j$, where $j$ is the source and $i$ the target. 
\end{ex}

\begin{defi}
\label{defpath}
Let $Q=(Q_V ,  Q_A , s , t)$ be a quiver. A path of $Q$ is a $p$-tuple of arrows $\pi = (a_1 , \dots , a_p) \in Q_A^p$ such that for every $j \in \{ 1 , \dots , p-1 \}$, $t(a_j) = s(a_{j+1})$. Its source (resp. target) $s(\pi)$ (resp. $t(\pi)$) is $s(a_1)$ (resp. $t(a_p)$) and its length is $p$.
\end{defi}

Every vertex of a quiver also defines a path of length zero. An oriented cycle is then a path of positive length  with same source and target. A quiver is said to be {\it acyclic} iff it contains no oriented cycle.

\begin{lemma}
\label{lemmaposet}
A quiver $Q$ without loop is acyclic iff the relation "$i \leq j$ iff $Q$ contains a path from $j$ to $i$" defines a partial order on its vertices.
\end{lemma}

\begin{proof}
The relation $\leq$ is reflexive thanks to the existence of paths of length zero and transitive by concatenation of paths. Now, if $i,j$ are vertices such that $i \leq j$ and $j \leq i$, we get by concatenation an oriented cycle containing $i$ and $j$. This forces $i=j$ under the hypothesis that $Q$ is acyclic. Conversely, since $Q$ contains no loop by hypothesis, any oriented cycle contains at least two vertices $i \neq j$. The latter satisfy $i \leq j$ and $j \leq i$, which contradicts the antisymmetry of the order.
\end{proof}

\begin{defi}
\label{defgrading}
A grading of a quiver $Q=(Q_V ,  Q_A , s , t)$ is an injection $\gr : Q_V \to \Z$ such that for every path $\pi$ of positive length, $\gr \circ s (\pi) > \gr \circ t (\pi) $. An arrow $a \in Q_A$ gets then graded by the difference $\gr \circ s (a) - \gr \circ t (a) $.
\end{defi}

\begin{lemma}
\label{lemmagrading}
A quiver is gradable iff it is acyclic.
\end{lemma}

\begin{proof}
Let $Q=(Q_V ,  Q_A , s , t)$ be a quiver. If it is gradable with grading $\gr : Q_V \to \Z$, then the relation "$i \leq j$ iff $\gr (i) \leq \gr (j)$" defines a total order on $Q_V$ which strictly decreases along paths of positive length, so that $Q$ cannot contain any oriented cycle. Conversely, if $Q$ is acyclic, then the relation "$i \leq j$ iff $Q$ contains a path from $j$ to $i$" defines a partial order on $Q_V$ by Lemma \ref{lemmaposet}. We then get a grading on $Q$ by sending any minimal element $i_0$ of $Q_V$ to zero and then by induction, for every $j \geq 0$, by sending any minimal element $i_{j+1}$ of $Q_V \setminus \{ i_0 , \dots , i_{j} \}$ to $j+1 \in \Z$.
\end{proof}

\begin{ex}
\label{exgrading}
Let $V$ be a finite subset of $\Z$. Then, the quiver $Q=(Q_V ,  Q_A , s , t)$ given by Example \ref{exquiver}, where $Q_A$ contains exactly one arrow for each pair $(i,j) \in Q_V^2$ such that $i \leq j$, is graded by the inclusion $V \hookrightarrow \Z$.
\end{ex}

\subsection{Shellability of a Morse tiling}
\label{subsecshellability}

Recall that a finite relative simplicial complex is the complement $S=K \setminus L$ for a subcomplex $L$ of a finite simplicial complex $K$, see \cite{Stan3,Stanbook,CKS}. We may assume that $L$ does not contain any maximal simplex of $K$, deleting them from $K$ otherwise, and call with some abuse $K = \overline{S}$ the underlying simplicial complex of $S$. Finite simplicial complexes themselves are special cases of relative ones for which $L$ is empty. Any Morse tile given by Definition \ref{defMorsetile} is another relative simplicial complex. We now recall the notion of Morse tilings on such finite relative simplicial complexes, as defined in \cite{SaWel2}, see also \cite{SaWel1,Wel1}. 

\begin{defi}[Definition $2.8$ of \cite{SaWel2} and $2.5$ of \cite{Wel1}]
\label{defMorsetiling}
A finite relative simplicial complex $S=K \setminus L$ is Morse tileable iff it admits a partition by Morse tiles such that for every $j \geq 0$, the union of tiles of dimension greater than $j$ is closed  in $S$. Such a partition is called a Morse tiling and the collection of simplices underlying these tiles, together with their faces, is the underlying simplicial complex $\overline{S} \subset K$ of $S$.
\end{defi}

Recall that a pure dimensional finite simplicial complex $K$ is classically said to be {\it shellable} whenever its maximal simplices  can be numbered $\sigma_1 , \dots , \sigma_N$ in a such way that for every $p \in \{1 , \dots , N-1 \}$, $\sigma_{p+1} \setminus (\sigma_1 \cup \dots \cup \sigma_p)$ is a non-empty union of codimension one faces of $\sigma_{p+1}$, see \cite{Bing,Bjo,Koz,Z}. This condition has been relaxed in several ways, including collapsibility \cite{Whi,AdipBen,AdipBen2}, constructibility \cite{Hach}, local constructibility \cite{BenZ,HachZ} or partitionability \cite{Stan3,Stanbook}. Any shelling or partitioning of a pure dimensional finite simplicial complex provides a Morse tiling, where the tiles $\sigma_{p+1} \setminus (\sigma_1 \cup \dots \cup \sigma_p)$ are basic, see Example \ref{exshelling} for further examples. Likewise, a special class of Morse tilings consists of the following shellable ones and our purpose in this section is to characterize when a tiling is shellable, proving Theorem \ref{theointro1}.

\begin{defi}[Definition $2.14$ of \cite{SaWel2} and $2.10$ of \cite{Wel1}]
\label{defMorseshelling}
A Morse tiled finite relative simplicial complex $S$ is shellable iff it admits a filtration $\emptyset = S_0 \subset S_1 \subset \dots \subset S_N = S$ by closed subsets such that for every $j \in \{1 , \dots , N\}$, $S_j \setminus S_{j-1}$ consists of a single tile. 
\end{defi}

When the tiling uses only basic tiles, this notion of shelling recovers, in fact slightly extends, the classical one, see Theorem $2.15$ and Remark $2.16$ of \cite{SaWel2}.

\begin{figure}[h]
   \begin{center}
    \includegraphics[scale=1]{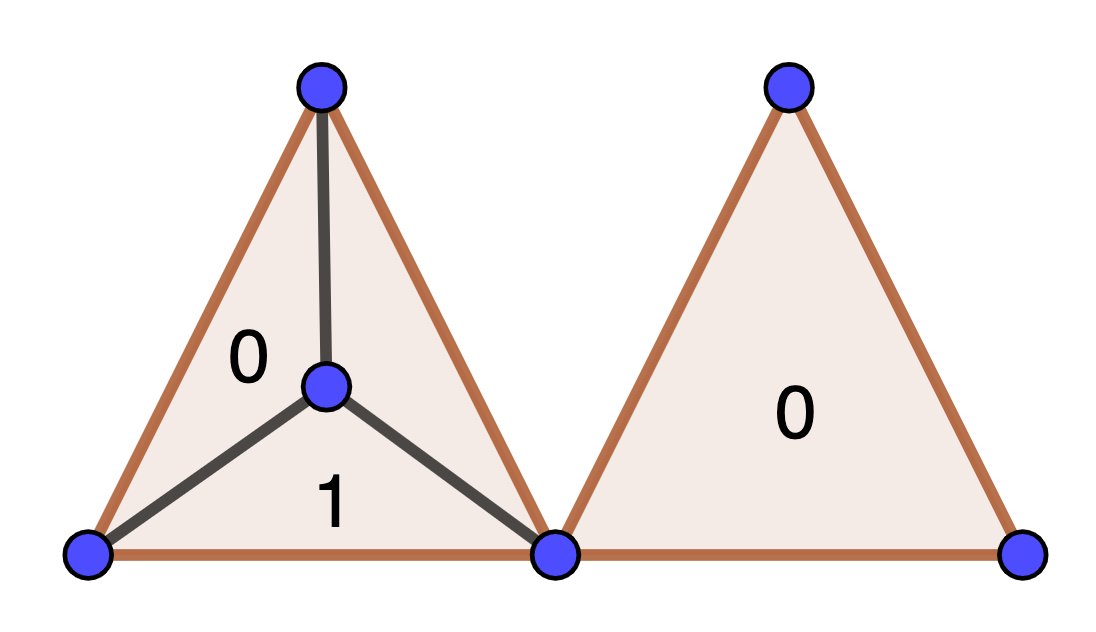}
    \caption{A Morse shelled complex in dimension two.}
    \label{figstellar}
      \end{center}
 \end{figure}

\begin{ex}
\label{exshelling} 
1) The boundary of every convex simplicial polytope is shellable by \cite{BruMan}, so that it is Morse shellable, the shelling using only basic tiles. 

2) The union of two simplices sharing a common face of codimension greater than one in both simplices is Morse shellable and any Morse shelling uses a closed simplex together with a regular Morse tile of order zero.

3) By Theorem $1.3$ of  \cite{SaWel2}, every closed triangulated surface is Morse shellable. 

4) By Theorem $1.3$ of \cite{Wel2}, every finite simplicial complex becomes Morse shellable after finitely many stellar subdivisions at maximal simplices, or also after a single barycentric subdivision. Figure \ref{figstellar} shows an example of such a Morse shelling taken out from \cite{Wel2}, where the numbers refer to the indices of its critical tiles. 
\end{ex}

More examples of Morse shellings can be found in \cite{Wel1}. We now observe that a Morse tiling supports a quiver which encodes the obstruction to shell it.

\begin{defi}
\label{defMorsequiver} 
The quiver of a Morse tiled finite relative simplicial complex $S$ has the tiles of $S$ as vertices and contains an arrow $T \to T'$ iff $\overline{T} \cap T' \neq \emptyset$. The arrows are labelled by the dimension of the face $\overline{T} \cap {T}'$ and the vertices by the order of the corresponding tiles. 
\end{defi}

In particular, the label of the target of an arrow cannot exceed the label of the arrow itself by more than one, see \S \ref{sectiles}.

\begin{figure}[h]
   \begin{center}
    \includegraphics[scale=1]{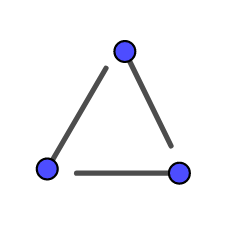}
    \caption{A non shellable tiling of $\partial \Delta_{2}$.}
    \label{figexotic}
      \end{center}
 \end{figure}

\begin{ex}
\label{exMorsequiver}
1) If $S = T_0 \sqcup \dots \sqcup T_{n+1}$ is the shelling of $\partial \Delta_{n+1}$ using one $n$-dimensional basic tile of each order, where $\Delta_{n+1}$ denotes the standard $(n+1)$-simplex, then its quiver is the comparability graph with vertices $T_0 , \dots , T_{n+1}$  and arrows $T_j \to T_i$ for every $0 \leq i \leq j \leq n+1$, see Corollary $4.2$ of \cite{SaWel1}.

2) If $S = T_1 \sqcup T_1 \sqcup T_1$  is the non-shellable tiling of $\partial \Delta_{2}$ using three basic tiles of order and dimension one, then its quiver is the oriented boundary of a triangle, see Figure \ref{figexotic}.
\end{ex}

We now prove Theorem \ref{theointro1}.

\begin{proof}[Proof of Theorem \ref{theointro1}]
Let us first assume that a Morse tiled finite relative simplicial complex $S$ is shelled and denote by $T_1 , \dots , T_N$ the shelling order on its tiles. Let $a : T_j \to T_i$ be an arrow of its quiver, so that $\overline{T}_j \cap T_i \neq \emptyset$. Since by definition the union $\cup_{k\leq j} T_k$ is closed in $S$, this forces $i \leq j$ so that the shelling order of $S$ induces a grading of its quiver by Definition \ref{defgrading}. We deduce from Lemma \ref{lemmagrading} that the latter is acyclic. 

Conversely, if the quiver of $S$ is acyclic, it is gradable by Lemma \ref{lemmagrading} and any grading provides a total order on the tiles of $S$. We may then label them $T_1 , \dots , T_N$ in increasing order. Then, for every $k \in \{ 1, \dots , N \}$ and every $1 \leq i \leq k < j \leq N$, $\overline{T}_i \cap T_j = \emptyset$, for otherwise the quiver of $S$ would contain an arrow $T_i \to T_j$ by Definition \ref{defMorsequiver} which is impossible by Definition \ref{defgrading}. We deduce that $S_k = \cup_{i\leq k} T_i$ is closed in $S$, so that the filtration $(S_k)_{k \in \{ 1, \dots , N \}}$ is a Morse shelling of $S$.
\end{proof}

\begin{rem}
\label{remDMT}
1) The acyclicity condition in Theorem \ref{theointro1} may be compared to the criterium given in Theorem $9.3$ of \cite{For}, up to which a discrete vector field on a simplicial complex is the gradient vector field of a discrete Morse function iff it contains no non-stationary closed path. Indeed, a Morse tiling on a finite simplicial complex encodes a class of compatible discrete vector fields and in the case of a shellable one, these are gradient vector fields of discrete Morse functions whose critical points are in one-to-one correspondance with the critical tiles, preserving the indices, see \cite{SaWel2} and \S \ref{secMorse}.

2) Every product of a sphere with a torus of positive dimension carries Morse tiled triangulations using only basic tiles, see  \cite{Wel1}. These have been obtained by subdividing the product of the shelled boundary of a simplex with several copies of the tiling given by Figure \ref{figexotic}.
These tilings cannot be shelled, so that their quivers all contain oriented cycles, as in the second part of Example \ref{exMorsequiver}.
\end{rem}

We also deduce from Theorem \ref{theointro1} criteria of partial shellability of a Morse tiling. 

\begin{defi}
\label{defqquiver} 
The $q$-quiver of a Morse tiled finite relative simplicial complex $S$ has the tiles of order $\leq q+1$ of $S$ as vertices and contains an arrow $T \to T'$ iff $\overline{T} \cap T'  \cap S^{(q)} \neq \emptyset$. 
\end{defi}

The $q$-skeleton of a Morse tiled finite relative simplicial complex $S=K \setminus L$ is the finite relative simplicial complex $S^{(q)} = K^{(q)} \setminus L^{(q)}$, that is the intersection with $S$ of the $q$-skeleton of its underlying simplicial complex $\overline{S} \subset K$. 

\begin{theo}
\label{theopartialshel}
The $q$-skeleton $S^{(q)}$ of a Morse tiled finite relative simplicial complex $S$ with acyclic $q$-quiver is Morse shellable and carries a filtration $\emptyset = S_0 \subset S_1 \cap S^{(q)} \subset \dots \subset S_N \cap S^{(q)} = S^{(q)}$ by closed subset such that for every $j \in \{1 , \dots , N\}$, $S_j \setminus S_{j-1}$ consists of a single tile of order $\leq q+1$ of $S$.
\end{theo}

Any filtration given by Theorem \ref{theopartialshel} is called a {\it partial shelling}. Let us first recall that skeletons of Morse tiles are Morse shellable, see Lemma $2.17$ and Theorem $2.18$ of \cite{SaWel2}.

 \begin{lemma}[Lemma $2.17$ of \cite{SaWel2}]
 \label{Lemmaskel} 
  Skeletons of Morse tiles are Morse shellable.
  \end{lemma}
  
  \begin{proof} 
 Let us first prove the result for codimension one skeletons of Morse tiles. Let $T = \sigma \setminus (\sigma_0 \cup \dots \cup \sigma_{k-1} \cup \mu)$ be an $n$-dimensional Morse tile of order $k$, where $\sigma$ denotes an $n$-simplex, $\sigma_0 , \dots , \sigma_n$ its codimension one faces and $\mu$ some higher codimensional face which contains the vertices opposite to $\sigma_0 , \dots , \sigma_{k-1}$ and which can be assumed not to contain the vertex opposite to $ \sigma_{k}$. Then, $\partial \sigma$ gets shelled by the tiles $T_0 \sqcup \dots  \sqcup T_n$, where $T_0 = \sigma_0 $ and for every $i \in \{ 1, \dots , n \}$, $T_i = \sigma_i \setminus (\sigma_0 \cup \dots \cup \sigma_{i-1})$. Likewise, the $(n-1)$-skeleton of $T$ gets shelled by the Morse tiles
$(T_k \setminus \mu) \sqcup T_{k+1} \sqcup \dots  \sqcup T_n$, proving the result in this case. Now, replacing $T$ with its $(n-1)$-skeleton, we deduce that its codimension two skeleton is Morse shellable, a shelling being obtained by concatenation of the shellings of the codimension one skeletons of each of the tiles $T_k \setminus \mu , T_{k+1} ,  \dots , T_n$. The general case then follows by finite induction.
\end{proof}

\begin{proof}[Proof of Theorem \ref{theopartialshel}]
Let $S$ be a Morse tiled finite relative simplicial complex with acyclic $q$-quiver $Q$ and let $\widetilde{S}$ be the union of tiles of order $\leq q+1$ of $S$. Then, $\widetilde{S}$ contains the $q$-skeleton of $S$ since a tile of order greater than $q+1$ does not contain any open face of dimension $\leq q$, see \S \ref{sectiles}. By Definition \ref{defqquiver}, the vertices of $Q$ are the tiles of $\widetilde{S}$ and it contains an arrow $T \to T'$ iff $\overline{T} \cap T' \cap S^{(q)}$ is non-empty. By Lemma \ref{lemmagrading}, it is gradable and any grading provides a total order on the tiles of $\widetilde{S}$, so that we may label them $T_1 , \dots , T_N$ in increasing order. Then, for every $1 \leq i \leq k < j \leq N$, $\overline{T}_i \cap T_j  \cap S^{(q)}= \emptyset$, for otherwise $Q$ would contain an arrow $T_i \to T_j$ which is impossible by Definition \ref{defgrading}. We deduce that ${S}^{(q)}_k = \big( \cup_{i\leq k} T_i \big) \cap S^{(q)}$ is closed in $S^{(q)}$, so that it provides the desired filtration by closed subsets of $S^{(q)}$. Now, by Lemma \ref{Lemmaskel}, the $q$-skeleton of every Morse tile is Morse shellable. Choosing such a Morse shelling on each $T_i \cap S^{(q)}$, $i \in \{ 1, \dots , N \}$, we get the result by concatenation of these shelling orders.
\end{proof}

\section{Spectral sequences}
\label{secspectral}

We now observe that any Morse shelling on a finite relative simplicial complex provides a way to compute its relative (co)homology via two spectral sequences.

\begin{theo}
\label{theospectral}
Let $S = K \setminus L$ be a Morse tiled finite relative simplicial complex having acyclic $q$-quiver. Then, any partial shelling given by Theorem \ref{theopartialshel} induces a spectral sequence which converges to the relative homology $H_* (K^{(q)},L^{(q)};G)$ (resp. relative cohomology $H^* (K^{(q)},L^{(q)};G)$) of its  $q$-skeleton with coefficients in any Abelian group $G$ and whose first page in gradings less than $q$ is the relative homology (resp. cohomology) of its critical tiles.
\end{theo}

By Proposition \ref{proprelhom}, only the critical tiles of indices less than $q$ contribute to the spectral sequences given by Theorem \ref{theospectral} in gradings less than $q$.

\begin{rem}
\label{remhomology}
The inclusion $i$ of the  pair of $q$-skeletons $(K^{(q)},L^{(q)})$ into the ambient pair of simplicial complexes $(K,L)$ induces an homomorphism $i_* : H_* (K^{(q)},L^{(q)} ; G) \to H_* (K,L ; G) $ (resp. $i^* : H^* (K,L ; G)  \to H^* (K^{(q)},L^{(q)} ; G)$) in homology (resp. cohomology) which is surjective (resp. injective) in degree $q$ and bijective in degrees less than $q$, see \cite{Munk}, so that the spectral sequences given by Theorem \ref{theospectral} converge to the relative  (co)homology of the relative simplicial complex in degrees less than $q$.
\end{rem}

\begin{proof}
We may assume that the finite simplicial complex $K$ coincides with $\overline{S} \subset K$, so that the subcomplex $L$ does not contain any maximal simplex of $K$.
Then, by Theorem \ref{theopartialshel}, the $q$-skeleton of $S$ inherits from any partial shelling a filtration $\emptyset = S_0 \subset S_1 \cap S^{(q)} \subset \dots \subset S_N  \cap S^{(q)}= S^{(q)}$ by closed subsets such that for every $p \in \{ 1, \dots , N \}$, $S_p \setminus S_{p-1} = T_p \cap S^{(q)}$ where $T_p$ is a tile of order $\leq q+1$ of $S$. This filtration induces a filtration $0 \subset C_* (K_1,L_1 , G) \subset \dots \subset C_* (K_N,L_N , G) = C_* (K^{(q)},L^{(q)} ; G)$ of simplicial chain complexes, where for every $p \in \{ 1, \dots , N \}$, $K_p = \overline{S}_p \cap K^{(q)}$ and $L_p = \overline{S}_p \cap  L^{(q)}$, and thus a spectral sequence which converges to the relative homology of the pair $(K^{(q)},L^{(q)})$, see chapter $XV$ of \cite{CarEil} for instance. The zero-th page of this spectral sequence is by definition
$E^0_{p, *-p} = C_* (K_p , K_{p-1} \cup L_p)$, $p \in \{ 1, \dots , N \}$, where we omit from now on to mention the group of coefficients $G$. Its first page is $E^1_{p, *-p} = H_* (K_p , K_{p-1}\cup L_p)$, $p \in \{ 1, \dots , N \}$ and its limit term is the relative homology of $(K^{(q)},L^{(q)})$ under the form $E^\infty_{p, *-p} = H_* (K^{(q)},L^{(q)})_p / H_*(K^{(q)},L^{(q)})_{p-1}$, $p \in \{ 1, \dots , N \}$, where $H_* (K^{(q)},L^{(q)})_p$ denotes the image of the inclusion homomorphism $H_* (K_p,L_p) \to H_* (K^{(q)},L^{(q)})$. These images indeed provide the induced filtration $0  \subset H_* (K^{(q)},L^{(q)})_1  \subset  \dots \subset H_* (K^{(q)},L^{(q)})_N = H_* (K^{(q)},L^{(q)})$ of the total relative homology of the $q$-skeleton. Moreover, the filtration being of finite length $N$, the spectral sequence degenerates at the $N$-th page, so that $E^\infty = E^N$. Now, for every $p \in \{ 1, \dots , N \}$, we deduce by excision an isomorphism $H_* (K_p , K_{p-1}\cup L_p) \cong H_* (\overline{T}_p \cap K^{(q)} , (\overline{T}_p \setminus T_p) \cap K^{(q)})$, see \cite{Munk}, and then by composition with the inclusion homomorphism, an isomorphism $H_* (K_p , K_{p-1}\cup L_p) \cong H_* (\overline{T}_p , \overline{T}_p \setminus T_p)$ as long as $* < q$, see Remark \ref{remhomology}. The homology part of the statement then follows from Proposition \ref{proprelhom}.

The proof of the cohomology part is similar and obtained from the filtration $0 \subset C^* (K_N , K_{N-1} \cup L_N ; G) \subset \dots \subset C^* (K_N, K_0 \cup L_N ; G) = C^* (K^{(q)} ,L^{(q)} ; G)$ of cochain complexes. It induces a spectral sequence which converges to the relative cohomology of the pair $(K^{(q)},L^{(q)})$ and whose zero-th page is by definition
$E_0^{p, *-p} = C^* (K_N , K_p \cup L_N) / C^* (K_N , K_{p+1} \cup L_N) = C^* (K_{p+1} , K_p \cup L_{p+1})$, $p \in \{ 0, \dots , N-1 \}$, where the last isomorphism is given by restriction of the cochains to $K_{p+1}$ and provides the short exact sequence of the triad $(K_N, K_{p+1} \cup L_N , K_p \cup L_N)$, and where we again omit to mention the coefficient group $G$ from now on. Its first page is thus $E_1^{p, *-p} = H^* (K_{p+1} , K_p \cup L_{p+1})$, $p \in \{ 0, \dots , N-1 \}$ and its limit term is the relative cohomology of $(K^{(q)},L^{(q)})$ under the form $E_\infty^{p, *-p} = H^* (K^{(q)},L^{(q)})_p / H^*(K^{(q)},L^{(q)})_{p+1}$, $p \in \{ 0, \dots , N-1 \}$, where $H^* (K^{(q)},L^{(q)})_p$ denotes the image of the inclusion homomorphism $H^* (K^{(q)} , K_p \cup L^{(q)}) \to H^* (K^{(q)}, L^{(q)})$. Now as before, we deduce by excision and Remark \ref{remhomology}, for every 
$p \in \{ 0, \dots , N-1 \}$, an isomorphism $H^* (K_{p+1} , K_p\cup L_{p+1}) \cong H^* (\overline{T}_{p+1} , \overline{T}_{p+1} \setminus T_{p+1})$ as long as $* < q$, so that the result follows from Proposition \ref{proprelhom}.
\end{proof} 

We now deduce.

\begin{proof}[Proof of Theorem \ref{theointro2}]
Theorem \ref{theointro2} follows from Theorem \ref{theospectral} combined with Theorem \ref{theointro1} and Proposition \ref{proprelhom}.
\end{proof}

\begin{ex}
\label{exspectral}
1) Any shelling of the boundary of a convex simplicial polytope given by \cite{BruMan} provides a Morse shelling which uses only basic tiles among which a single closed simplex and a single open one. The spectral sequences given by Theorem \ref{theointro2} degenerate at the first page and compute the (co)homology of a sphere.

2) The Morse shellings of any closed connected triangulated surface given by Theorem $1.3$ of \cite{SaWel2} can also be chosen to use a single closed simplex by construction, see Theorem $4.1$ of \cite{Wel2}. Then, the spectral sequences given by Theorem \ref{theointro2} again degenerate at the first page whenever the Morse shelling uses a single open simplex as well and they compute then the (co)homology of the surface as the number of critical tiles of each index it uses.

3) In the case of the Morse shelled octahedron depicted in Figure \ref{figOctahedron} -taken out from \cite{Wel1}-, where the numbers appearing in Figure \ref{figOctahedron} indicate the shelling order, the spectral sequences given by Theorem \ref{theointro2} degenerate at the second page. Indeed, the (co)homology classes of the pair of critical tiles with shelling orders $6$ and $7$, which have consecutive indices one and two, get killed already at this second page.
\begin{figure}[h]
   \begin{center}
    \includegraphics[scale=1]{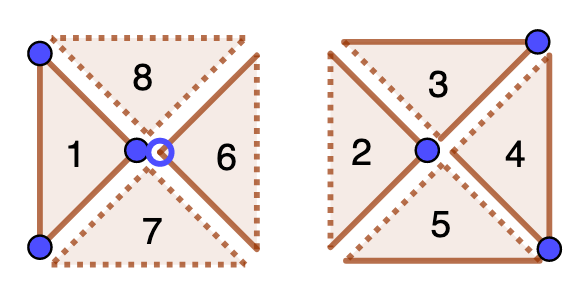}
    \caption{A Morse shelling on the octahedron.}
    \label{figOctahedron}
      \end{center}
 \end{figure}
 \end{ex}
 
 \begin{rem}
\label{remfinale}
1) It follows from Theorem \ref{theospectral} that the $k$-th Betti number of a Morse tiled simplicial complex with acyclic $q$-quiver, $q>k$, is bounded from above by the number of critical tiles of index $k$ of the tiling. In the case of a Morse shellable complex, this recovers Corollary $1.5$ of  \cite{SaWel2}, which has been deduced from discrete Morse theory.

2) Theorem \ref{theointro2} provides a way to compute the (co)homology of a Morse shellable simplicial complex using much smaller (co)chain complexes than the simplicial ones, as the discrete Morse complexes introduced by R. Forman would do, see \cite{For}. For instance, performing a large number of barycentric subdivisions on such a simplicial complex increases by a lot the dimensions of its simplicial (co)chain complexes, whereas by Theorem $1.1$ of  \cite{SaWel2}, see also \cite{Wel2}, these inherit Morse shellings using a constant number of critical tiles.
\end{rem}

\section{Relation with discrete Morse theory}
\label{secMorse}

We already proved in \cite{SaWel2} that every Morse tiling on a finite simplicial complex encodes a class of compatible discrete vector field which, in the case of a shellable one, are all gradient vector fields of discrete Morse functions. Let us now add that such a discrete Morse function can be chosen in such a way that the filtration given by the Morse shelling, which generates the spectral sequences, coincides with a filtration by sublevel sets of the function. Recall that a subcomplex $L$ of a finite relative simplicial complex $K$ is an {\it elementary collapse} of $K$ if and only if $K \setminus L$ consists of a maximal relative simplex of $K$ together with one of its codimension one face which is not contained in any other maximal simplex of $K$. It is a {\it collapse} of $K$ iff it is obtained from $K$ after finitely many elementary collapses, see for instance \cite{AdipBen,Whi}. 

\begin{prop}
\label{propcollapse}
A critical Morse tile of index $k$ collapses onto any of its $k$-dimensional face, whereas a regular Morse tile collapses onto the empty set. 
\end{prop}
In the case of a critical tile of maximal index, no collapse is needed at all in Proposition \ref{propcollapse}.

\begin{proof}
We proceed by induction on the dimension $n$ of the tiles. If $n=1$, a closed one-simplex collapses onto any of its vertices and a closed simplex deprived of one vertex collapses onto the empty set.
Let $T = \sigma \setminus (\sigma_0 \cup \dots \cup \sigma_{k-1} \cup \mu)$ be an $n$-dimensional Morse tile of order $k$, where $\sigma$ denotes an $n$-simplex, $\sigma_0 , \dots , \sigma_n$ its codimension one faces and $\mu$ some higher codimensional face. If $T$ is critical, let $\theta$ be any of the remaining $k$-dimensional faces of $\sigma$ ; it contains $\mu$.  Then, $\mu$
contains the vertices opposite to $\sigma_0 , \dots , \sigma_{k-1}$ and $\theta$ can be assumed not to contain the vertex opposite to $ \sigma_{k}$. We deduce that $\partial \sigma$ gets shelled by the tiles $T_0 \sqcup \dots  \sqcup T_n$, where $T_0 = \sigma_0 $ and for every $i \in \{ 1, \dots , n \}$, $T_i = \sigma_i \setminus (\sigma_0 \cup \dots \cup \sigma_{i-1})$. If $T$ is of order $n+1$, it is an open simplex and no collapse is needed to get Proposition \ref{propcollapse}. Otherwise, depriving $T$ of its open maximal face together with $T_n$, we deduce that it collapses onto $(T_k \setminus \mu) \sqcup T_{k+1} \sqcup \dots  \sqcup T_{n-1}$. By the induction hypothesis, it collapses then onto
$T_k \setminus \mu$. Since $T$ is critical, $\dim (\mu ) = k-1$ so that $T_k \setminus \mu$ is critical as well and contains $\theta$ by hypothesis. It thus collapses onto $\theta$ by the induction hypothesis. If $T$ is regular, $\mu$ is empty and we deduce likewise that $T$ collapses onto the empty set. Hence the result.
\end{proof}

In the framework of discrete Morse theory \cite{For}, a sequence of elementary collapses given by  Proposition \ref{propcollapse} is encoded by some discrete vector field $V$ such that for every collapsed pair of faces $ \sigma < \tau$ of the tile $T$, $V(\sigma) = \tau$ and such that $V$ vanishes on the remaining $k$-dimensional face in case $T$ is critical of index $k$. We say that a discrete vector field on a Morse tiled finite simplicial complex is {\it compatible} with the tiling if and only if it restricts to each tile of the tiling to such a discrete vector field given by  Proposition \ref{propcollapse}, compare \cite{SaWel2}.

\begin{theo}
\label{theoMorse}
Let $S=K \setminus L$ be a Morse shelled finite relative simplicial complex and $\emptyset = S_0 \subset S_1 \subset \dots \subset S_N = S$ be the corresponding filtration, where for every $j \in \{ 1, \dots , N \}$, $S_j \setminus S_{j-1}$ is a single tile of the shelling. Then, for every discrete vector field $V$ compatible with the shelling, there exists a discrete Morse function $f$ on $K= \overline{S}$ whose gradient vector field coincides with $V$ on $S$ and such that  for every $j \in \{ 0, \dots , N \}$, $f^{-1} (]- \infty , j])=S_j \cup L$.
\end{theo}

\begin{proof}
Let $T_1, \dots , T_N$ be the Morse tiles of the shelling and $V$ be any discrete vector field on $S$ compatible with the shelling. For every $j \in \{ 1, \dots , N \}$, let us consider a sequence of elementary collapses given by Proposition \ref{propcollapse} and encoded by the restriction of $V$ to $T_j$. We denote by  $(\tau^j_1 , \sigma^j_1), \dots , (\tau^j_{k_j} , \sigma^j_{k_j})$ the pairs of collapsed faces of $T_j$ in the order of collapses, so that $(\tau^j_1 , \sigma^j_1)$ is collapsed first and $(\tau^j_{k_j} , \sigma^j_{k_j})$ collapsed at the final step. Let us moreover assume that $\tau^j_{k_j} = \sigma^j_{k_j}$ in the case $T_j$ is critical, so that $V(\tau^j_{k_j} )=0$ in this case, whereas $V (\sigma^j_i) = \tau^j_i$ in all the other ones. Let the restriction of $f$ to $L$ be any discrete Morse function taking non positive values. We then set, for every $j \in \{ 1, \dots , N \}$ and every $i \in \{ 1, \dots , k_j \}$, $f(\sigma^j_i) = f(\tau^j_i) = j-1+\frac{1}{i}$. By construction, it is a discrete Morse function with gradient vector field $V$ on $S$ and such that for every $j \in \{ 0, \dots , N \}$, $f^{-1} (]- \infty , j])=S_j \cup L$. Hence the result.
\end{proof}

In the case of a finite simplicial complex, the critical faces of the Morse functions given by Theorem \ref{theoMorse} are in one-to-one correspondence with the critical tiles of the shelling, preserving the index, compare \cite{SaWel2}. The spectral sequences given by Theorem \ref{theointro2} are thus associated to the filtration by sublevel sets of some discrete Morse functions on the simplicial complex. They can be computed using discrete Morse theory in the following way. Any discrete Morse function $f$ given by Theorem \ref{theoMorse} induces a cellular decomposition of $K$ using one cell for each critical face of $f$, see \S $7$ of \cite{For}. Each page of the spectral sequence can then be computed using this decomposition. Indeed, each cell of the decomposition provides a (co)cycle representative of the relative (co)homology of the corresponding critical tile and thus a representative of the corresponding (co)homology class in the first page. It remains a (co)cycle representative of this (co)homology class as long as the differential of the corresponding page maps it to a shelling degree where the boundary of the cell contains a simplex, which must then be a critical point of degree one less of the Morse function. The latter sits inside a critical tile of the shelling. Then, whether or not this homology class survives to higher pages depends on whether or not this boundary component can be annihilated by another cell, which is encoded by the discrete Morse complex.

It is however not needed to choose a compatible discrete Morse function to compute the spectral sequences, as the case of the octahedron given in Example \ref{exspectral} shows. They can be computed directly from the shelling.

\begin{rem}
By Theorem $8.10$ of \cite{For}, the differentials of the Morse complex of a discrete Morse function can also be computed as a finite sum of gradient paths between critical points of consecutive indices. If the Morse function is chosen to be compatible with a Morse shelling, then any such gradient path provides an oriented path in the quiver of the shelling, whose source and target are critical tiles of consecutive indices.
\end{rem}

\bibliographystyle{abbrv}

Univ Lyon, Universit\'e Claude Bernard Lyon 1, CNRS UMR 5208, Institut Camille Jordan, 43 blvd. du 11 novembre 1918, F-69622 Villeurbanne cedex, France

{\tt welschinger@math.univ-lyon1.fr.}
\end{document}